\documentclass[a4paper,12pt]{amsart}

\usepackage{fullpage,amsxtra,amssymb,amsthm,amsmath}
\usepackage[latin1]{inputenc}
\numberwithin{equation}{section}

\renewcommand{\leq}{\leqslant}
\renewcommand{\geq}{\geqslant}

\def\stacksum#1#2{{\stackrel{{\scriptstyle #1}}
{{\scriptstyle #2}}}}

\newcommand{\Cc}{\mathbf{C}}

\newcommand{\Zz}{\mathbf{Z}}

\newcommand{\Rr}{\mathbf{R}}

\newcommand{\Qq}{\mathbf{Q}}
\newcommand{\Fp}{\mathbf{F}}

\newcommand{\mods}[1]{\,(\mathrm{mod}\,{#1})}

\newcommand{\ideal}[1]{\mathfrak{{#1}}}

\newcommand{\ra}{\rightarrow}

\DeclareMathOperator{\ord}{ord}

\DeclareMathOperator{\Tr}{Tr}

\DeclareMathOperator{\End}{End}

\DeclareMathSymbol{\gena}{\mathord}{letters}{"3C}
\DeclareMathSymbol{\genb}{\mathord}{letters}{"3E}

\theoremstyle{plain}
\newtheorem{theorem}{Theorem}[section]

\newtheorem{lemma}[theorem]{Lemma}
\newtheorem{corollary}[theorem]{Corollary}

\newtheorem{proposition}[theorem]{Proposition}

\theoremstyle{remark}
\newtheorem{remark}[theorem]{Remark}

\theoremstyle{definition}

\newcommand{\wg}{w}
\newcommand{\pr}{\mathfrak{p}}

\begin{document}

\title{Weil numbers generated by other Weil numbers and torsion fields
  of abelian varieties}
\author{E. Kowalski}
\address{Universit\'e Bordeaux I - A2X\\
351, cours de la Lib\'eration\\
33405 Talence Cedex\\
France}
\email{emmanuel.kowalski@math.u-bordeaux1.fr}
\begin{abstract}
Using properties of the Frobenius eigenvalues, we show that, in a
precise sense, ``most'' isomorphism classes of (principally polarized)
simple abelian varieties over a finite field are characterized, up to
isogeny, by the sequence of their division fields, and a
similar result for ``most'' isogeny classes. Some global cases are
also treated. 
\end{abstract}

\maketitle

\section{Introduction}

Let $q=p^k$ be a power of 
a prime number and let $\pi$ be a $q$-Weil
number, i.e., an algebraic integer such that for every automorphism
$\sigma$ of $\Cc$ we have $|\sigma(\pi)|=\sqrt{q}$. Let $\Phi_{\pi}$ be the
multiplicative group generated inside $\Cc^{\times}$ by the Galois-conjugates
of $\pi$. We are interested in the multiplicative group structure of
$\Phi_{\pi}$ and particularly in the set (say $\wg_{\pi}$) of $q$-Weil numbers
inside $\Phi_{\pi}$. Clearly, $\wg_{\pi}$ contains all the conjugates of
$\pi$, and we wish to know when there is equality.
\par
Our motivation relates to abelian varieties over finite
fields. Let $A/\Fp_q$ be such an abelian variety over a field with $q$
elements. Weil proved that all eigenvalues of the Frobenius
endomorphism $\pi_A$ of $A$ are $q$-Weil numbers. We denote by
$\Phi_A$ the multiplicative group that they generate. If $A$ is
simple, then $\pi_A$ ``is'' an algebraic integer and
$\Phi_A=\Phi_{\pi_A}$. We also denote by $\wg_A$ the set of $q$-Weil
numbers in $\Phi_A$.
\par
In~\cite[Th. 3.4]{kow}, it is shown that given abelian
varieties $A/\Fp_q$ and $B/\Fp_q$, the condition
$\Phi_{A}=\Phi_{B}$ is equivalent with the fact that, for all $n$
coprime with some integer $N$ (which may depend on $A$ and $B$), we
have $\Fp_q(A[n])=\Fp_q(B[n])$, in which case we say that $A$ and $B$
are isokummerian. 
In the case of a simple variety,
therefore, if $\wg_{\pi_A}$ is 
reduced to the conjugates of $\pi_A$, any variety $B$ satisfying the
above condition must be isogenous to a power of $A$.
\par
Some precise results (e.g., for product of elliptic curves) are given
in~\cite{kow}, and also an example
due to Serre of two abelian varieties $A$ and $B$ over a finite field,
simple and non-isogenous over 
$\bar{\Fp}_q$, such that $\Phi_A=\Phi_B$. (In particular, in such a case
$\wg_{A}$ contains strictly the set of conjugates
of $\pi_A$). 
\par
Our main result shows that for ``most'' abelian varieties
over finite fields, $\wg_{A}$ is the
set of eigenvalues of $\pi_A$. What ``most'' means has to be
specified, of course, and there are actually at least two natural ways
of doing this, taking $A$ up to isomorphism or isogeny. We will
consider both possibilities, using a common lemma and specialized
arguments. In the
isomorphism case, we use rather deep results of Mumford and Chavdarov
on what ``most'' isomorphism classes of abelian varieties over $\Fp_q$
look like; it is quite appealing that we use here both $p$-adic
methods having to do with ordinariness and $\ell$-adic methods
related to monodromy of $\ell$-adic sheaves. Those also allow us to
derive some results for abelian varieties over number fields, although
they are conditional on ordinarity assumptions. In the isogeny case,
the method is more elementary, based on lattice-point counting, using
results of Howe and DiPippo, and the multidimensional large sieve
inequality.
\par
\medskip
We now state precisely our main result. Let $g\geq 1$ be an integer,
$q=p^k$ with $p$ prime and $k\geq 1$. We introduce the following:
\begin{itemize}
\item $A_g(q)$ is the set of isomorphism classes of principally
polarized abelian varieties of dimension $g$ defined over $\Fp_{q}$;
\item $I_g(q)\subset A_g(q)$ is the subset of those varieties $A$ such
  that any $B$ isokummerian with $A$ is isogenous to a power of $A$;
\item $\mathcal{A}_g(q)$ is the set of isogeny classes of abelian
  varieties of dimension $g$ defined over $\Fp_{q}$;
\item $\mathcal{I}_g(q)\subset \mathcal{A}_g(q)$ is the subset of
  isogeny classes such that any $B$ isokummerian to $A$ is (isogenous
  to) a power of $A$. 
\end{itemize}

\begin{theorem}\label{th-main}
Let $g\geq 1$ and $q=p^k$ with $p$ prime and $k\geq 1$. 
We have
\begin{equation}\label{eq-isom}
\lim_{n\ra +\infty}\frac{|I_g(q^n)|}{|A_g(q^n)|}=1
\end{equation}
and
\begin{equation}\label{eq-iso-g}
\lim_{n\ra +\infty}\frac{|\mathcal{I}_g(q^n)|}{|\mathcal{A}_g(q^n)|}
=1.
\end{equation}
\end{theorem}

\textbf{Acknowledgments.} The results of Chavdarov~\cite{chavdarov}
which are crucial for this paper were mentioned by N. Katz during a
lecture; shortly afterward a question by U. Zannier made me realize
that those results could be quite useful to study $\Phi_A$ and
$\wg_A$ and improve on~\cite{kow}. I thank them both for these lucky
coincidences...

\section{Determination of $\wg_{\pi}$ in a special case}

In this section we consider only $q$-Weil numbers, and give a
criterion for $\wg_{\pi}$ to be reduced to the conjugates of $\pi$.
For simplicity we assume that $\pi$ does not have real conjugates,
hence $\pi$ is of even degree $2g$. Let $K_{\pi}\subset \bar{\Qq}$ be
the Galois closure of $\Qq(\pi)$. For every conjugate $\pi_i$ of
$\pi$, $q/\pi_i$ is also a conjugate of $\pi$; if we fix an embedding
$\bar{\Qq}\subset \Cc$, we have $q/\pi_i=\bar{\pi}_i$, the complex
conjugate of $\pi_i$.

\begin{proposition}\label{pr-2}
Let  $p$ be prime, $q=p^k$ with $k\geq 1$. Let $\pi$ be a $q$-Weil
number such that
$[\Qq(\pi):\Qq]=2g$. Let $G$ denote the 
Galois group of the Galois closure $K_{\pi}$ of $\Qq(\pi)$. Let
$(\pi_i,\bar{\pi}_i)$, $1\leq i\leq g$, 
be the Galois conjugates of $\pi$ in $K_{\pi}$, in complex conjugate
pairs. Assume that:
\par
\emph{(1)} For all $i$, $1\leq i\leq g$, $\pi_i$ and $\bar{\pi}_i$ are
coprime in $K_{\pi}$.
\par
\emph{(2)} For all $i$, $1\leq i\leq g$, there exists
$\sigma_i\in G$ such that $\sigma_i(\pi_i)=\pi_i$ and
$\sigma_i(\pi_j)=\bar{\pi}_j$ for $j\not=i$. 
\par
Then $\wg_{\pi}$ is the set of conjugates of $\pi$.
\end{proposition}


\begin{proof}
Let $c\in G$ be the restriction of complex conjugation, so that
$c(\pi_i)=\bar{\pi}_i$ and $c(\bar{\pi}_i)=\pi_i$ for every
$i$. Thus setting
$$
\sigma_{i,j}=c\sigma_i\sigma_j\in G
$$
we have for $i\not=j$ the relations
\begin{equation}\label{eq-sigmaij}
\sigma_{i,j}(\pi_i)=\pi_i,\quad \sigma_{i,j}(\pi_j)=\pi_j,\quad
\sigma_{i,j}(\pi_k)=\bar{\pi}_k\text{ for } k\notin \{i,j\}.
\end{equation}
\par
Fix a prime ideal $\pr$ in $K_{\pi}$ dividing $(p)$. Since $\pr\mid p\mid
q=\pi_i\bar{\pi}_i$ and $\pi_i$, $\bar{\pi}_i$ are coprime by
assumption, we see that $\pr$ divides one and only one of $\pi_i$ and
$\bar{\pi}_i$. We 
renumber/pair the conjugates so that $\pr\mid\pi_i$, and $\pr\nmid
\bar{\pi}_i$  for $1\leq i\leq g$. Notice this doesn't affect the
existence of $\sigma_i$, $\sigma_{i,j}$ with properties as stated for
the new numbering.
\par
Let
$\nu=v_{\pr}(q)\geq 1$ where $v_{\pr}$ is the valuation on $K_{\pi}$
associated to $\pr$. Notice that by coprimality again we have
\begin{equation}\label{eq-nu}
\nu=v_{\pr}(q)=v_{\bar{\pr}}(q)=v_{\pr}(\pi_i\bar{\pi}_i)=v_{\pr}(\pi_i)
=v_{\bar{\pr}}(\bar{\pi}_i).
\end{equation}
\par
Let now $\alpha\in \Phi_{\pi}$ be a $q$-Weil number. We can write
$$
\alpha=q^m\prod_{1\leq i\leq g}{\pi_i^{n_i}},
$$
with $n_i\in\Zz$. We deduce from $\alpha\bar{\alpha}=q$ that
\begin{equation}\label{eq-p2}
2m+n_1+\cdots+n_g=1,
\end{equation}
and from this we notice in particular that the sum $n_1+\cdots+n_g$ can not be
zero, in particular not all the $n_i$ can be zero.
\par
We have $v_{\pr}(\alpha)\geq 0$, $v_{\bar{\pr}}(\alpha)\geq 0$, which
translate to
$$
\nu(m+n_1+\cdots+n_g)\geq 0,\quad \nu m\geq 0.
$$
\par
Dividing by $\nu\geq 1$, summing and comparing with~(\ref{eq-p2}), we
see that one of $m$ and 
$m+n_1+\cdots+n_g$ is equal to $0$ and the other is equal to $1$.
\par
Now we consider $\alpha\alpha^{\sigma_i}$. This is an algebraic
integer and therefore $v_{\pr}(\alpha\alpha^{\sigma_i})\geq 0$, 
$v_{\bar{\pr}}(\alpha\alpha^{\sigma_i})\geq 0$. We have
$$
\alpha\alpha^{\sigma_i}=q^{2m+n_1+\cdots+\cdots+n_g-n_i}\pi_i^{2n_i}=
q^{1-n_i}\pi_i^{2n_i},
$$
so using~(\ref{eq-nu}) these two conditions translate to
\begin{align*}
v_{\pr}(\alpha\alpha^{\sigma_i})=\nu(1+n_i)\geq 0\\
v_{\bar{\pr}}(\alpha\alpha^{\sigma_i})=\nu(1-n_i)\geq 0,
\end{align*}
which means that $n_i\in \{0,-1,1\}$ for all $i$.
\par
Now consider $\alpha\alpha^{\sigma_{i,j}}$ with $i\not=j$. We have
$$
\alpha\alpha^{\sigma_{i,j}}=q^{1-n_i-n_j}\pi_i^{2n_i}\pi_j^{2n_j},
$$
by~(\ref{eq-sigmaij}), hence the integrality conditions
$v_{\pr}(\alpha\alpha^{\sigma_{i,j}})\geq 0$, 
$v_{\bar{\pr}}(\alpha\alpha^{\sigma_{i,j}})\geq 0$ mean
\begin{align*}
v_{\pr}(\alpha\alpha^{\sigma_{i,j}})=\nu(1+n_i+n_j)\geq 0\\
v_{\bar{\pr}}(\alpha\alpha^{\sigma_{i,j}})=\nu(1-n_i-n_j)\geq 0.
\end{align*}
\par
The first of these shows that at most one $n_i$ can be equal to $-1$; the
second that at most one $n_j$ can be equal to $1$. Both of these can
not occur because that would give $n_1+\cdots +n_g=n_i+n_j=1-1=0$,
which is impossible. So either there exists exactly one $i$ with
$n_i=1$, and the other $n_j$ are $0$, which gives $\alpha=\pi_i$
(because one must have $m=0$, 
$m+n_1+\cdots+n_g=1$); or there exists exactly one $j$ with
$n_j=-1$ (and the other $n_i$ are $0$), which gives
$\alpha=\bar{\pi}_j$ (because then $m=1$, $m+n_1+\cdots+n_g=0$). 
\end{proof}

\begin{remark}
Since we actually solved the equations in terms of
the parameters $(m,n_i)$ uniquely (for a given $\alpha$), we have also
proved that $(q, \pi_i)$, $1\leq i\leq q$, form a free generating set
of $\Phi_{\pi}$ under the assumptions of the proposition. In
particular, the rank of $\Phi_{\pi}$ is then equal to $g+1$.
\end{remark}

\begin{remark}
Proposition~\ref{pr-2} also applies to prove that if
$A=E_1\times \cdots\times E_k$ is a product of pairwise geometrically
non-isogenous elliptic curves over a finite field $\Fp_q$ with $q$
elements, the only $q$-Weil numbers in $\Phi_A$ are the conjugates of
the Frobenius elements for the $E_i$ (see~\cite[Th. 3.4,
(5)]{kow}). It also gives back in this case the lemma of Spiess used
to prove this statement in loc. cit.
\end{remark}

To apply Proposition~\ref{pr-2} to a simple abelian variety $A/\Fp_q$
with Frobenius $\pi_A$, we need criteria for the two conditions
involved. Here we start by Condition~(1), which has to do with the
``behavior at $p$'' (since all primes dividing $\pi$ are above $p$ in
$K_{\pi}$) of the Frobenius of $A$.
\par
Recall that an abelian variety $A/k$ of
dimension $g$ over a field $k$ of characteristic $p$ is called
ordinary if $|A[p](\bar{k})|=p^g$, which is the maximal number of
$p$-torsion points there can be in characteristic $p$. If $k=\Fp_q$ is
a finite field with $q$ elements, then $A$ is ordinary if and only if
the middle coefficient of the characteristic polynomial of Frobenius
is coprime with $q$ (see e.g.~\cite{dipippo-howe}). The
following lemma is certainly well-known. 

\begin{lemma}\label{lm-ord}
Let $q=p^k$ with $p$ prime, $k\geq 1$, let $A/\Fp_q$ be a simple
ordinary abelian variety. Then for any eigenvalue $\pi$ of the
Frobenius of $A$, we have $(\pi,q/\pi)=1$ in the Galois closure of
$\Qq(\pi)$. 
\end{lemma}

\begin{proof}
Let $(\pi_i,q/\pi_i)$, $1\leq i\leq g$, be the conjugates of
$\pi$ with $\pi_1=\pi$. Assume there exists $k$ and a prime ideal
$\pr\mid p\mid q$ in $K_{\pi}$ dividing both $\pi_k$ and
$q/\pi_k$. The middle coefficient $b$ 
of the characteristic polynomial of Frobenius is given by
$$
b=\sum_{\stacksum{I,J\subset \{1,\ldots, g\}}{|I|+|J|=g}}{\prod_{i\in
    I}{\pi_i}\prod_{j\in J}{\frac{q}{\pi_j}}}. 
$$
In this sum, if $k\in I$ or $k\in J$, we have 
$$
\pr \mid \prod_{i\in I}{\pi_i}\prod_{j\in J}{\frac{q}{\pi_j}}
$$
by assumption. Otherwise, $I$ and $J$ are both chosen inside
the set $\{1,2,\ldots, g\}-\{k\}$ with $g-1$ elements, and
$|I|+|J|=g$. Thus $I\cap J\not=\emptyset$. If $i\in I\cap J$, then 
$$
\pr\mid q=\pi_i\cdot  q/\pi_i\mid \prod_{i\in I}{\pi_i}
\prod_{j\in J}{\frac{q}{\pi_j}}.
$$
Thus we find that $\pr$ divides all the terms in the sum giving
$b$. Since $b\in\Zz$, this means $p\mid b$. Hence $A$ is not ordinary,
and the result follows by contraposition.
\end{proof}

This implies that any ordinary abelian variety satisfies Condition~(1)
of Proposition~\ref{pr-2}.  Note that the examples of Serre
in~\cite{kow} are not ordinary (since their endomorphism rings are not
commutative, which is another consequence of ordinarity, see
e.g.~\cite[\S 7]{waterhouse}).  
\par
In analogy with $A_g(q)$, $\mathcal{A}_g(q)$, we now denote
\begin{itemize}
\item $A_{g}^{\ord}(q)\subset A_{g}(q)$ the set of isomorphism classes of
  principally polarized ordinary abelian varieties of dimension $g$
  defined over $\Fp_{q}$; 
\item $\mathcal{A}_{g}^{\ord}(q)$ the set of isogeny classes of ordinary
abelian varieties of dimension $g$ defined over $\Fp_{q}$.
\end{itemize}
\par
Now we come to Condition~(2), where there is also a simple sufficiency
criterion.

\begin{lemma}\label{lm-w2g}
Let $\pi$ be a $q$-Weil number such that $[\Qq(\pi):\Qq]=2g$ and such
that the Galois group of $K_{\pi}$ over $\Qq$ is isomorphic to
$W_{2g}$, the Weyl group of $Sp(2g)$. Then $\pi$ satisfies
Condition~(2) of Proposition~\ref{pr-2}.
\end{lemma}

\begin{proof}
Recall that $W_{2g}$ (the Galois group of a ``generic'' polynomial $P$
of degree $2d$ such that $X^{2d}P(1/X)=P(X)$) can be identified with
the group of permutations of $g$ pairs $(2i-1,2i)$, $1\leq i\leq g$,
such that the couples $\{2i-1,2i\}$ are stable. In the case of the
Galois group of $K_{\pi}$, the pairs can be identified with the pairs
of conjugates $(\pi_i,q/\pi_i)$, which shows that $G$ can be
identified with a subgroup of $W_{2g}$. If it is equal to $W_{2g}$,
the existence of the required elements $\sigma_i$ is obvious. 
\end{proof}

\begin{corollary}\label{cor-crit}
Let $q=p^k$ with $p$ prime and $k\geq 1$. For any simple ordinary abelian
variety $A/\Fp_q$ of dimension $g$ such that the Galois group $G$ of
$K_{\pi_A}$ is isomorphic to $W_{2g}$, the set of $q$-Weil numbers in
$\Phi_A$ is equal to the set of conjugates of $\pi_A$.
\end{corollary}

This is immediate from Proposition~\ref{pr-2}, Lemma~\ref{lm-ord}
and Lemma~\ref{lm-w2g}. Again we denote:
\begin{itemize}
\item $B_{g}(q)\subset A_{g}(q)$ the set of isomorphism classes of
  absolutely simple principally polarized abelian varieties of dimension $g$
  defined over $\Fp_{q}$ such that the Galois group of $K_{\pi_A}$ is
  isomorphic to $W_{2g}$; 
\item $\mathcal{B}_{g}(q)\subset \mathcal{A}_g(q)$ the set of isogeny
  classes of absolutely simple abelian varieties of dimension $g$
  defined over $\Fp_{q}$ such that the Galois group of $K_{\pi_A}$ is
  isomorphic to $W_{2g}$.
\end{itemize}

\section{General abelian varieties up to isomorphism}

We now apply Proposition~\ref{pr-2} to ``generic'' isomorphism classes
of abelian varieties of dimension $g$. More precisely, one has to
consider (for instance) the moduli space $A_g$ of abelian varieties of
dimension $g$ with a principal polarization, which is known to be
irreducible of dimension $g(g+1)/2$ over $\Zz$.
\par
For Condition~(1) of Proposition~\ref{pr-2}, we use
Lemma~\ref{lm-ord}. It is known that generic abelian varieties are
ordinary (see~\cite{oort-norman}\footnote{~The result is 
  attributed to Mumford~\cite{mumford}, although the author confesses
  that he doesn't see
  that statement in this paper of Mumford.}.) Thus in
$A_g$, there exists a dense Zariski open subset $U\subset A_{g}$ such that
the polarized abelian variety parameterized by any $u\in U$ is
ordinary. (See also~\cite[\S 5]{chambert-loir} for a sketch; roughly
speaking, ordinarity is an open condition, and we know that ordinary
abelian varieties of any dimension exist, for instance products of
ordinary elliptic curves).

\begin{proposition}\label{pr-ord-isom}
Let $q=p^k$ with $p$ prime, $k\geq 1$. We have
$$
\lim_{n\ra +\infty}{\frac{|A_{g}^{\ord}(q^n)|}{|A_g(q^n)|}}=1.
$$
\end{proposition}

\begin{proof}
Mumford's result gives this for the corresponding counting of
isomorphism classes of principally polarized abelian varieties with
some rigidifying structure; then one deduces the statement above by
dividing by the number of choices for the rigidifying data, and
dealing with possible extra automorphisms, as done for instance
in~\cite[10.7,11.3]{katz-sarnak}. 
\end{proof}

\begin{remark}
This is much weaker than what the result of Mumford implies: since the
space of abelian varieties is of dimension $g(g+1)/2$ and the space of
non-ordinary abelian varieties must be of dimension $\leq 
g(g+1)/2-1$, we have for $n\geq 1$
\begin{align*}
|A_g(q^n)|&=q^{ng(g+1)/2}+O(q^{n(g(g+1)/2-1)}),\\
|A^{\ord}_g(q^n)|&=|A_g(q^n)|+O(q^{n(g(g+1)/2-1)}).
\end{align*}
\end{remark}

Condition~(2) is not so easy to treat. We use Lemma~\ref{lm-w2g}, and
the crucial fact is that Chavdarov~\cite{chavdarov} has shown that 
``most'' abelian varieties $A/\Fp_{q^n}$ with $n\ra +\infty$ are
simple and satisfy the assumptions of that lemma. 

\begin{proposition}\label{pr-ch}
Let $q=p^k$ with $p\not=2$ prime and $k\geq 1$. Then we have
$$
\lim_{n\ra +\infty}{\frac{|B_g(q^n)|}{|A_g(q^n)|}}=1.
$$
\end{proposition}

\begin{proof}
This follows from~\cite[Th. 2.1]{chavdarov}, applied to a suitably
``rigidified'' universal family of principally polarized abelian
varieties of dimension $g$ over $\Fp_q$, after eliminating as before
the extra factor counting the rigidifying parameters (compare
again~\cite[11.3]{katz-sarnak}). The monodromy groups modulo $\ell$
involved in applying Chavdarov's Theorem are as large as possible for
$\ell>2$ because (for instance), it is already the case for the
families of jacobians of hyperelliptic curves considered
in~\cite[Th. 11.0.4]{katz-sarnak}, which are the same as those
in~\cite[Ex. 2.4]{chavdarov} (this is where characteristic $\not=2$
enters). This is due to J.K. Yu (unpublished). See also below for more
discussion of these examples. 
\end{proof}

We now deduce from Proposition~\ref{pr-ord-isom} and
Proposition~\ref{pr-ch} the first main result of this paper.

\begin{theorem}\label{th-iso}
Let $q=p^k$ with $p$ prime and $k\geq 1$. For $n\geq 1$, let
$C_g(q^n)$ be set of isomorphism classes of principally polarized
absolutely simple abelian varieties $A/\Fp_{q^n}$ of dimension $g$
such that $\Phi_A\simeq \Zz^{g+1}$ and $\wg_{A}$ is equal to the set of
conjugates of $\pi_A$. Then we have 
$$
\lim_{n\ra +\infty}{\frac{|C_g(q^n)|}{|A_g(q^n)|}}=1.
$$
\end{theorem}

Informally: ``most'' abelian varieties $A$ of dimension $g$ over
$\Fp_{q^n}$ with $n$ large are simple, ordinary, the group $\Phi_A$ is
isomorphic to $\Zz^{g+1}$ and the only $q^n$-Weil numbers in $\Phi_A$
are $\pi_A$ and its conjugates.
\par
By the criterion stated in the introduction for two varieties to be
isokummerian over a finite field, we see that this theorem is
equivalent with the first part~(\ref{eq-isom}) of
Theorem~\ref{th-main}. 


\begin{remark}
As in~\cite{chavdarov} or~\cite{katz-sarnak}, it would be very
interesting to have a corresponding result with $n=1$ and $q=p\ra
+\infty$; and (as in those cases) this seems very hard. 
\par
On the other hand, introducing some analytic ideas (a bilinear form
estimate for representations of $\Fp_{\ell}$-adic sheaves and
``old-fashioned'' large sieve as in~\cite{gallagher} and 
Section~\ref{sec-iso}),
it is possible to improve Proposition~\ref{pr-ch} in some cases (in
particular, if $g$ satisfies $p>2g+1$) to obtain a sharper
estimate 
$$
|I_g(q^n)|=q^{ng(g+1)/2}+O(q^{n(g(g+1)/2-\gamma)}(\log q^n))
$$
for $\gamma=(10g^2+6g+8)^{-1}$; see~\cite[Cor. 6.4]{sieve-monodromy}. 
\end{remark}

If one does not wish to deal with the moduli space, one can apply
Chavdarov's theorem to any algebraic family of principally polarized
abelian varieties over 
a finite field $\Fp_q$ for which the monodromy group mod $\ell$ is
equal to $Sp(2g,\Zz/\ell\Zz)$ for almost all $\ell$, provided one can check
that ordinarity is generic in that family. The simplest example are
provided by taking an algebraic family of curves and then the
associated jacobian family, which has a canonical principal
polarization. If one takes the universal family of curves, then the
generic ordinarity is a result of Miller, who gives explicit examples
of ordinary curves of every genus and characteristic, so that the result
follows from the openness of ordinarity and the irreducibility of the
moduli space of curves. The fact that the corresponding monodromy
group is $Sp(2g)$ follows again in characteristic $\not=2$ from the
examples of families of hyperelliptic curves
of~\cite[Ex. 2.4]{chavdarov} (see also~\cite[10.2]{katz-sarnak}).
\par
It is natural to want to give similar explicit equations of families
of curves which are both generically ordinary and have monodromy
$Sp(2g)$. However note that Miller's families
$$
\begin{cases}
y^2=x^{2g+1}+tx^{g+1}+x\text{ if } p\nmid g,\\
y^2=x^{2g+2}+tx^{g+1}+x\text{ if } p\mid g
\end{cases}
$$
fail the monodromy test (because they fail the diophantine
irreducibility test, see~\cite[Lemma 10.1.15]{katz-sarnak}, as a
simple computation shows). On the other hand, the author couldn't find
references to the ordinarity for the families with large monodromy of
loc. cit. For the moment, we merely state the following fairly easy
result: 
\begin{proposition}
Let $p\geq 3$ be a prime number, $g\geq 2$ an integer. Put $\delta=1$
if $p\mid g$, $\delta=0$ otherwise.
\par
\emph{(1)} The $2$-parameter family $T$ of smooth projective curves of
genus $g$ over $\Fp_p$ given by compactification of the affine family 
$$
T_{t,u}\,:\,y^2=(x-u)(x^{2g+\delta}+tx^g+1)
$$
over the open subset
$$
U=\{(t,u)\in\mathbf{A}^2\,\mid\, u^{2g+\delta}+tu^g+1\not=0\}\subset 
\mathbf{A}^2/\Fp_p
$$
is generically ordinary and 
has geometric monodromy group modulo $\ell$ equal to $Sp(2g,\Fp_{\ell})$
for $\ell>2$, $\ell\not=p$.
\par
\emph{(2)} In particular, there exists $k\geq 1$ and $t_0\in
\Fp_{p^k}$ such that the $1$-parameter family $S$ of curves of genus $g$
over $\Fp_{p^k}$ given by
$$
S_u\,:\, y^2=(x-u)(x^{2g+\delta}+t_0x^g+1)
$$
with $u\in U_{t_0}=\{(t_0,u)\in U\}\subset \mathbf{A}^1/\Fp_{p^k}$ is
generically ordinary and has 
geometric monodromy group modulo $\ell$ equal to $Sp(2g,\Fp_{\ell})$
for $\ell>2$, $\ell\not=p$.
\par
\emph{(3)} If $p\nmid g$ and $g$ is even, one can in fact take $t_0=0$,
so the family 
$$
S_u\,:\, y^2=(x-u)(x^{2g}+1)
$$
with $u\in U_0=\mathbf{A}^1-\mu_{2g}$, where $\mu_{2g}$ is the group
of $2g$-roots of unity, is generically ordinary and has
geometric monodromy group modulo $\ell$ equal to $Sp(2g,\Fp_{\ell})$
for $\ell>2$, $\ell\not=p$.
\end{proposition}

\begin{proof}
Note that for $u=0$, the family $T$ specializes to Miller's family,
and therefore the open set of ordinarity for $T$ is not empty, hence
dense. Moreover, for
any fixed $t_0\not=\pm 1$, the family $S_u=T_{u,t_0}$ is of the form
$$
y^2=f_{t_0}(x)(x-u)
$$
with $f_{t_0}$ a monic polynomial of degree $2g$ which has distinct
roots in $\bar{\Fp}_p$, as a simple computation shows. Hence it is of
the form considered in~\cite[Ex. 2.4]{chavdarov}
and~\cite[10.1]{katz-sarnak}, and therefore has 
the required monodromy. As the monodromy groups can only become
smaller by taking such a $1$-parameter subfamily of $T$, the result
follows. 
\par
Now generic ordinarity for $T$ implies that for some
$t_0\in\mathbf{A}^1(\bar{\Fp}_p)-\{\pm 1\}$ at least the restricted
subfamily $S_u=T_{u,t_0}$ with $u$ as parameter must contain an
ordinary curve. As it still has the required monodromy groups, the
existence result~(2) follows.
\par
When $p\nmid g$, and $t=t_0=0$, Miller's curve has equation
$$
y^2=x^{2g+1}+x.
$$
The recipe in~\cite[\S 2]{miller} for computing the (dual of the)
Hasse-Witt matrix for this curve shows that it is invertible, hence
the curve is ordinary, if for every $u$, $0\leq u\leq
g-1$, there exist unique integers $r$, $t\geq 0$ such that
\begin{gather*}
r+t=\frac{p-1}{2}\\
2gt+\frac{p+1}{2}+u=p(v+1),
\end{gather*}
where $v$ is uniquely determined by $0\leq v\leq g-1$ and the congruence
\begin{equation}\label{eq-1}
\frac{p+1}{2}+u\equiv p(v+1)\mods{g}.
\end{equation}
It is easy to see (following Miller's argument) that the equations for
$r$ and $t$ have at most one solution, and that this solution exists if
and only if 
$$
v+1\equiv u+\frac{p+1}{2}\mods{2}.
$$
If $g$ is even, then~(\ref{eq-1}) implies this.\footnote{\ On the
  other hand, if $g$ is odd, the residue modulo $2$ of this expression
will take both values.}
\end{proof}

We now come to some global consequences that follow also from other
results of Chavdarov's paper. Those have the virtue of concerning
individual abelian varieties, as the exceptions become a set of primes
of density $0$ which does not affect (for instance) the Isogeny
Theorem. 

\begin{proposition}
Let $F/\Qq$ be a number field. Let $g\geq 1$ be $2$, $6$, or an odd
integer. Let $A/F$ be an abelian variety of dimension $g$ such that
$\End(A)=\Zz$ and such that the set of primes of good reduction
$\ideal{p}$ of $F$ where the reduction of $A$ modulo $\ideal{p}$ is
ordinary is of density $1$. Then for any abelian variety  $B/F$, $B$
is isokummerian to $A$ if and only if $B$ is isogenous to a power of
$A$. 
\end{proposition}

\begin{proof}
By Chavdarov's ``horizontal'' version of his result
(\cite[Cor. 6.9]{chavdarov}), the assumption of $A$ ensures that for all prime
ideals $\ideal{p}$ in a set of primes of density $1$, the reduced
variety $A_{\ideal{p}}/\Fp_{\ideal{p}}$ is ordinary, absolutely simple
and its Frobenius has Galois group $W_{2g}$. By
Corollary~\ref{cor-crit}, it follows that $B_{\ideal{p}}$ must be
isogenous to a power of $A_{\ideal{p}}$ for any such $\ideal{p}$. The
dimension of $B$ fixes a $k\geq 1$ such that $B_{\ideal{p}}\simeq
A_{\ideal{p}}^k$ for all primes in a set of density $1$. Then by
Faltings's Isogeny Theorem, it follows that $B\simeq A^k$ over $F$. 
\end{proof}

The assumption of ordinarity at almost all places for varieties with
$\End(A)=\Zz$ is widely expected to hold, but few results are
known. For elliptic curves, it is quite easy, but this case of the
proposition is already treated in~\cite{kow} without this
assumption. Here is another situation that can be treated
unconditionally (compare with Ogus's theorem quoted
in~\cite[Th. 6.3]{chambert-loir}):

\begin{proposition}
Let $A/\Qq$ be an abelian surface over $\Qq$ with $\End(A)=\Zz$. Then
the set of primes of good ordinary reduction for $A$ is of density
$1$. Hence any $B/\Qq$ is isokummerian to $A$ if and only if $B$ is
isogenous to a power of $A$.
\end{proposition}

\begin{proof}
We use Serre's $\ell$-adic methods~\cite{serre-chebo}. Let $\ell$ be a
prime and $\rho_{\ell}\,:\, G_{\Qq}\ra Sp(4,\Qq_{\ell})$ the
$\ell$-adic representation associated to $A$. Serre has shown (this is
already used in the proof of Chavdarov's horizontal theorem) that the
image of $\rho_{\ell}$ is dense. Consider the exterior
square $\sigma_{\ell}=\wedge^2\rho_{\ell}$. It is an $\ell$-adic
representation of rank $6$ and ``weight'' $1$, and it is faithful, so
the closure $G_{\ell}$ of the image of 
$\sigma_{\ell}$ is again isomorphic to $Sp(4,\Qq_{\ell})$. Moreover,
for any prime 
$p\not=\ell$ of good reduction, the properties of $\rho_{\ell}$ and
standard algebra show that the trace of the image by $\sigma_{\ell}$ of a
Frobenius element $\sigma_p$ at $p$ is the middle coefficient $b_2$ of
the characteristic polynomial of the Frobenius of $A$ modulo $p$. Hence,
by the characterization of ordinarity already stated, $A$ has ordinary
reduction at $p$ if and only if this trace
$\Tr\sigma_{\ell}(\sigma_p)$ is not divisible by $p$. 
\par
However, by the Riemann Hypothesis for $A$ modulo $p$, we have
$$
|\Tr\sigma_{\ell}(\sigma_p)|\leq 6p,
$$
so if $A$ is not ordinary at $p$, the trace must belong to
set $\{-6p,-5p,\ldots, 0,p,\ldots, 6p\}$. Let $t$ be any
of these thirteen values. We 
claim that the set of primes $p$ for which
$\Tr\sigma_{\ell}(\sigma_p)=t$ is of density $0$. Clearly this implies
the proposition.
\par
The proof of the claim is easy: since
$\det\sigma_{\ell}(\sigma_p)=p^4$, if $p$ satisfies the stated
condition then we have
$$
\sigma_p\in X_t=\{g\in G_{\ell}\,\mid\, 
(\Tr g)^4-t^4\det g=0\}.
$$
Using $G_{\ell}\simeq Sp(4,\Qq_{\ell})$ and simple computations, it is
easy to see that $X_t$ is a closed subset of $G_{\ell}$ of Minkowski
dimension $<\dim Sp(4)$ (see~\cite[\S 3]{serre-chebo} for the definition
of Minkowski, or $M$-dimension). Hence by Theorem 10 of loc. cit., the
set of primes with $\sigma_p\in X_t$ is of density $0$.
\end{proof}



In a general higher dimensional situations (over $\Qq$, say), the
non-ordinary primes are such that the trace $t$ of the $g$-th exterior
power of the 
representation on the Tate module are divisible by $p$, which for
$g\geq 3$ allows an unbounded number of values of $t$ (for $g=3$, $\Tr
\wedge^3\rho_{\ell}(\sigma_p)=pk$ with $|k|\leq 20\sqrt{p}$). Even
using explicit forms of the Chebotarev density theorem (on GRH) to
detect each value, the uniformity is not sufficient to obtain any
non-trivial result.



\begin{remark}
In~\cite{kow}, the question of the ``splitting behavior'' of a simple
abelian variety $A/\Qq$ at all primes is also raised: is it true that the
reduction modulo $p$ of $A$ remains simple for almost all $p$? In
fact, the ``horizontal'' statements of Chavdarov already deal with
this. For instance, this holds if $A/\Qq$ has the property that the
Galois group of the field $\Qq(A[\ell])$ generated by the points of
$\ell$-torsion of $A$ is 
equal to $Sp(2g,\Zz/\ell\Zz)$ for $\ell$ large.
\end{remark}

\section{General abelian varieties up to isogeny}
\label{sec-iso}

Since Weil numbers, ordinarity, and having Galois group $W_{2g}$ are
all isogeny-invariant properties of abelian varieties, it is natural
to ask for analogs of the results of the previous section for
isogeny classes of abelian varieties, instead of isomorphism classes. 
Going directly from one to the other is not easy, since finding
the number of isomorphism classes in an isogeny class is a
quite delicate question, typically related with class numbers (as can be
seen most easily in the case of elliptic curves),
see~\cite[\S 4.3]{waterhouse}.
\par
However, we can use results of DiPippo and Howe to deal directly
with isogeny classes. Note then that it is not necessary to introduce 
a polarization. This is rather satisfactory since not all isogeny
classes contain a principally polarized one; see for
instance~\cite[Th. 1.3]{howe-ppol}; however it is proved there (Th. 1.2) that
any isogeny class of odd-dimensional abelian varieties over a finite field
contains a principally polarized one.

\begin{proposition}\label{pr-ord-iso}
Let $g\geq 1$, $q=p^k$ with $p$ prime and $k\geq 1$. We have
$$
\lim_{n\ra  +\infty}
{\frac{|\mathcal{A}^{\ord}(q^n)|}{|\mathcal{A}_g(q^n)|}}=1,
$$
and
\begin{equation}\label{eq-iso-space}
|\mathcal{A}_g(q^n)|,\ |\mathcal{A}_g^{\ord}(q^n)|\sim v_g
\frac{\varphi(q^n)}{q^n}q^{ng(g+1)/4}.
\end{equation}
for some constant $v_g>0$.
\end{proposition}

This is proved by DiPippo and Howe in~\cite{dipippo-howe} (see Theorem
1.1), in fact in a much more precise form. Note in particular that
this says intuitively that the ``dimension''
of the space of isogeny classes of abelian varieties of dimension $g$
is $g(g+1)/4$, half that of the moduli space.

\begin{proposition}\label{pr-w2g-isogeny}
We have
$$
\lim_{n\ra +\infty}{\frac{|\mathcal{B}_g(q^n)|}{|\mathcal{A}_g(q^n)|}}=1.
$$
\end{proposition}

Using Lemma~\ref{lm-w2g}, this shows that the analogue of
Theorem~\ref{th-iso} holds for isogeny classes, and therefore that the
second part~(\ref{eq-iso-g}) of Theorem~\ref{th-main} holds.
\par
To prove Proposition~\ref{pr-ord-iso}, DiPippo and Howe identify the
set of isogeny classes considered with a set of lattice points in a
region $V_{g,n}\subset \Rr^g$. We argue similarly for
Proposition~\ref{pr-w2g-isogeny}, except that we do not need to be so
precise because we only look for an upper bound on the number of
isogeny classes with ``smaller'' Galois group, which is
a question of probabilistic Galois theory. It is straightforward to
adapt here
the method of Gallagher~\cite{gallagher} based on the large sieve
inequality. It has already been shown, using those methods, that 
self-reciprocal polynomials of degree $2g$ and bounded height have
generically $W_{2g}$ as Galois group (see~\cite{dds}), but our
parameter set is different.
\par
Let $A/\Fp_q$ be an abelian variety of dimension $g$ over a number
field. The  characteristic polynomial of Frobenius $f_A$ of $A$ is of
degree $2g$ with real roots of even multiplicity and complex roots
arising in pairs $(\alpha,q/\alpha)$. Therefore one can write
$$
f_A=(X^{2g}+q^g)+a_1(X^{2g-1}+q^{g-1}X)+\cdots+a_gX^g,
$$
with $a_i\in\Zz$. 
To $A$ we associate the vector $a=(a_1,\ldots,
a_g)\in \Zz^g$.

\begin{lemma}\label{lm-rgq}
Let $g\geq 1$, let $q=p^k$ with $p$ prime and $k\geq 1$. For any
abelian variety $A/\Fp_q$, the vector $a$ above satisfies $a\in
\Zz^g\cap R_{g,q}$ where
$$
R_{g,q}=\Bigl\{(x_1,\ldots, x_g)\in \Rr^g\,\mid\,
|x_i|\leq \binom{g}{i}q^{i/2}\Bigr\}.
$$
\end{lemma}

\begin{proof}
This is obvious by the Riemann Hypothesis and the definition of $a_i$. 
\end{proof}

The analytic ingredient we need is the following consequence of the
large sieve inequality.

\begin{lemma}\label{lm-lg-sieve}
Let $g\geq 1$. For $1\leq i\leq g$, let $X_i\geq 1$ and let
$$
R=\{(x_1,\ldots, x_g)\in \Zz^g\,\mid\, |x_i|\leq X_i\text{ for }
1\leq i\leq g\}\subset \Rr^g.
$$
Let $y\geq 2$ and for all primes $p\leq y$, let $\Omega(p)\subset
(\Zz/p\Zz)^g$ be a finite set of cardinality $\omega(p)$. Let
$$
P(y)=\sum_{p\leq y}{\omega(p)p^{-g}}
$$
and for any $a\in \Zz^g$ let $P(a,y)$ denote the number of $p\leq y$
such that $a\mods p\in \Omega(p)$.
Then we have
$$
\sum_{a\in R}{(P(a,y)-P(y))^2}\ll P(y)\prod_{j=1}^k{(X_i+y^2)},
$$
the implied constant depending only on $g$.
\end{lemma}

\begin{proof}
We derive this from the following multidimensional (trigonometric)
large sieve inequality: for any finite set of vectors $Y\subset
(\Rr/\Zz)^g$ such that $\max \|\alpha_k-\beta_k\|\geq \delta$ for two
elements $\alpha\not=\beta$ in $Y$ (where $\|x-y\|$ is the distance in
$\Rr/\Zz$), and for any complex numbers $f(x)$ defined for $x\in R$, we have
\begin{equation}\label{eq-lg}
\sum_{\alpha\in Y}{\Bigl|\sum_{x\in R}{f(x)e(\langle x,\alpha\rangle)}
\Bigr|^2}\ll \prod_{k=1}^g{(X_i+\delta^{-1})}\sum_{x\in R}{|f(x)|^2},
\end{equation}
where the implied constant depends only on $g$. This is a special case
of~\cite[Th. 1]{huxley}.
\par
To obtain the lemma from this, proceed as in Lemma~A
of~\cite{gallagher}, which we repeat for convenience: let $\chi_p$ be
the characteristic function of $\Omega(p)$, and expand it in Fourier
series 
$$
\chi_p(a)=\sum_{\alpha\in(\Zz/p\Zz)^g}{\hat{\chi}_p(\alpha)
e(\langle a,\alpha\rangle/p)}\text{ with }
\hat{\chi}_p(\alpha)=p^{-g}\sum_{a\in \Omega(p)}
{e(\langle -a,\alpha\rangle/p)}.
$$
Thus we have
\begin{equation}\label{eq-positivity}
\hat{\chi}_p(0)=p^{-g}\omega(p),\quad
\sum_{\alpha\not=0}{|\hat{\chi}_p(\alpha)|^2}\leq
\sum_{\alpha}{|\hat{\chi}_p(\alpha)|^2}=p^{-g}\omega(p).
\end{equation}
We have for $a\in R$
\begin{equation}\label{eq-pay}
P(a,y)=\sum_{p\leq y}{\sum_{\alpha\in(\Zz/p\Zz)^g}{\hat{\chi}_p(\alpha)
e(\langle a,\alpha\rangle/p)}}=
P(y)+\sum_{p\leq y}{\sum_{\alpha\not=0}{\hat{\chi}_p(\alpha)
e(\langle a,\alpha\rangle/p)}}.
\end{equation}
Denote by $R(a,y)$ the inner sum. We now write by Cauchy's inequality
and~(\ref{eq-positivity}) 
\begin{align*}
\sum_{a\in R}{|R(a,y)|^2}&=
\sum_{p\leq y}{\sum_{\alpha\not=0}{\hat{\chi}_p(\alpha)
\sum_{a\in R}{R(a,y)e(\langle a,\alpha\rangle/p)}}}\\
&\leq
\Bigl(
\sum_{p\leq y}{\sum_{\alpha\not=0}{|\hat{\chi}_p(\alpha)|^2}}
\Bigr)^{1/2}
\Bigl(
\sum_{p\leq y}{\sum_{\alpha\not=0}{\Bigl|
\sum_{a\in R}{R(a,y)e(\langle a,\alpha\rangle/p)}
\Bigr|^2}}
\Bigr)^{1/2}\\
&\leq P(y)^{1/2}
\Bigl(
\sum_{p\leq y}{\sum_{\alpha\not=0}{\Bigl|
\sum_{a\in R}{R(a,y)e(\langle a,\alpha\rangle/p)}
\Bigr|^2}}
\Bigr)^{1/2}
\end{align*}
and applying the trigonometric large sieve inequality~(\ref{eq-lg})
with the trivial spacing estimate for distinct vectors $\alpha/p$,
$\beta/q\in (\Rr/\Zz)^g$, $p$, $q\leq y$, this gives 
$$
\sum_{a\in R}{|R(a,y)|^2}\ll 
P(y)^{1/2}\Bigl(
\prod_{k=1}^g{(X_i+y^2)}\sum_{a\in R}{|R(a,y)|^2}\Bigr)^{1/2},
$$
so
$$
\sum_{a\in R}{|R(a,y)|^2}\ll P(y)\prod_{k=1}^g{(X_i+y^2)}.
$$
As, by~(\ref{eq-pay}), we have
$$
\sum_{a\in R}{(P(a,y)-P(y))^2}=
\sum_{a\in R}{|R(a,y)|^2},
$$
we are done.
\end{proof}

\begin{proof}[Proof of Proposition~\ref{pr-w2g-isogeny}]
First, for any $g$-tuple $a=(a_1,\ldots, a_g)$ in a ring $R$, we denote
$$
h_a=X^g+a_{1}X^{g-1}+\cdots +a_{g-1}X+a_1\in R[X]
$$
and
$$
f_a=X^gh_a(X+X^{-1})\in R[X].
$$
\par
Let $A$ be an abelian variety and $f_A$ the characteristic polynomial
of Frobenius for $A$ and $G$ the Galois group of its splitting field,
which can be seen (in possibly many ways) as a subgroup of $W_{2g}$.
By Lemma~2 of~\cite{dds}, we have $G=W_{2g}\subset \mathfrak{S}_{2g}$
if $G$ contains a $2$-cycle, a $4$-cycle, a $(2g-2)$-cycle and a
$2g$-cycle.  
\par
For $\ell\in \{2,4,2g-2,2g\}$, let $E_{\ell}$ denote the
number of lattice points $a=(a_1,\ldots, a_g)$ in the region $R_{g,q}$
defined in Lemma~\ref{lm-rgq} such that the polynomial $f=f_a$ 
is either reducible or such that the Galois group $G_a$ of the splitting
field of $f$, seen as a subgroup of $W_{2g}$ again, does not contain
an $\ell$-cycle. By the observation above and
Lemma~\ref{lm-rgq}, it follows that the number $E$ of isogeny classes of
abelian varieties $A/\Fp_q$ with $f_A$ not having Galois group
$W_{2g}$ satisfies
$$
E \leq E_2+E_4+E_{2g-2}+E_{2g}.
$$
\par
For each $\ell$, we know from classical algebraic number theory (see
e.g.~\cite[\S 61]{vdw}) that if the polynomial $f_a$ reduces modulo some prime
$p$ to a polynomial $f_a\mods p\in \Fp_p[X]$  which factorizes as a
product of $2g-\ell$ distinct linear factors and a single irreducible
polynomial of degree $\ell$, then $G_a$ contains an
$\ell$-cycle. Therefore, choosing $y\geq 2$ arbitrary and putting 
\begin{multline*}
\Omega(p)=\{a=(a_1,\ldots, a_g)\in (\Zz/p\Zz)^g\,\mid\,
f_a\mods{p}\text{ factorizes as $2g-\ell$ distinct}\\
\text{ linear factors, and one irreducible factor of degree $\ell$} \}
\end{multline*}
for $p\leq y$, we see that for $a$ such that $G_a$ does not contain 
an $\ell$-cycle we have $f_a\mods{p}\notin \Omega(p)$ for all $p\leq
y$. With notation as in Lemma~\ref{lm-lg-sieve} with
$X_i=\binom{g}{i}q^{i/2}$ (so $R=R_{g,q}$), we have therefore
$P(a,y)=0$, and the 
large sieve inequality implies by positivity that
$$
E_{\ell} P(y)^2\ll P(y)\prod_{1\leq i\leq g}{(q^{i/2}+y^2)}
$$
where the implied constant depends only on $g$.
However by Lemma~3 of~\cite{dds}
(see p. 269, or compare~\cite[p. 96, l. 10]{gallagher}) we have for
$y\geq 3$ the lower bound
$$
P(y)=\frac{C_{\ell}}{|W_{2g}|}\pi(y)+O(\log\log y)
\gg \pi(y),
$$
where $C_{\ell}$ is the number of $\ell$-cycles in $W_{2g}$, where the
implied constant depend only on $g$. 
Thus we get by the Prime Number Theorem (Chebychev's elementary
lower-bound estimate suffices) that
$$
E_{\ell}\ll \prod_{1\leq i\leq g}{(q^{i/2}+y^2)}
y^{-1}\log y.
$$
We choose $y^2=q^{1/2}$, so that
$$
\prod_{1\leq i\leq g}{(q^{i/2}+y^2)}\leq 2^gq^{g(g+1)/4}
$$
and
$$
E_{\ell}\ll q^{g(g+1)/4-1/4} \log q,
$$
hence
$$
E\ll q^{g(g+1)/4-1/4} \log q
$$
with an implied constant depending only on $g$. By
comparison with~(\ref{eq-iso-space}), we see that
Proposition~\ref{pr-w2g-isogeny} is proved.
\end{proof}

\begin{remark}
The bound obtained from the large sieve estimate may seem quite poor
because of the choice of a rather small $y$, constrained by the
smallest $X_i$. One may certainly expect that having a small Galois
group would be of ``codimension'' at least $1$, which would mean
essentially $E\ll q^{g(g+1)/4-1/2}$. There is a similar discrepancy
between what is proved and what is expected in other problems of
probabilistic Galois theory.
\end{remark}

\begin{remark}
In contrast with the isomorphism case, the results above do not yield
examples of ``thinner'' families of isogeny classes which would be
ordinary and have the $W_{2g}$ as associated Galois group. Most
notably, it is by no means clear how to prove the analogue
of~(\ref{eq-iso-g}) where the isogeny classes are jacobians of curves
of genus $g$ (equivalently, where arbitrary Weil numbers are replaced
by those associated with curves). Distinguishing jacobians among
abelian varieties over a finite field is a deep unsolved problem.
\end{remark}


\begin{thebibliography}{CCCC}


\bibitem[CL]{chambert-loir}
A. Chambert-Loir: \textit{Cohomologie cristalline: un survol},
Exposition. Math.  16  (1998), 333--382. 

\bibitem[C]{chavdarov}
N. Chavdarov: \textit{The generic irreducibility of the numerator of
  the zeta function in a family of curves with large monodromy}, Duke
Math. J. 87 (1997), 151--180.

\bibitem[DDS]{dds}
S. Davis, W. Duke and X. Sun: \textit{Probabilistic Galois theory of
reciprocal polynomials}, Exposition. Math.  16  (1998),  no. 4,
263--270.

\bibitem[DH]{dipippo-howe}
S. DiPippo and E. Howe: \textit{Real polynomials with all roots on the
unit circle and abelian varieties over finite fields}, J. Number
Theory 73 (1998), 426--450; Corrig., J. Number Theory 83 (2000), 182.

\bibitem[G]{gallagher}
P.X. Gallagher: \textit{The large sieve and probabilistic Galois
  theory}, in Proc. Sympos. Pure Math., Vol. XXIV, Amer. Math. Soc.
(1973), 91--101.

\bibitem[H]{howe-ppol}
E. Howe: \textit{Principally polarized ordinary abelian varieties over
finite fields}, Trans. Amer. Math. Soc. 347 (1995), 2361--2401.

\bibitem[Hu]{huxley}
M.N. Huxley: \textit{The large sieve inequality for algebraic number
  fields}, Mathematika 15 (1968), 178--187.


\bibitem[KS]{katz-sarnak}
N. Katz and P. Sarnak: \textit{Random matrices, Frobenius eigenvalues
  and monodromy}, A.M.S Colloquium Publ. 45, 1999.

\bibitem[K1]{kow}
E. Kowalski: \textit{Some local-global applications of Kummer theory},
manuscripta math. 111 (2003), 105--139.

\bibitem[K2]{sieve-monodromy}
E. Kowalski: \textit{The large sieve, monodromy and zeta functions of
  curves}, preprint (2005), \verb|arXiv:math/NT0503714|.

\bibitem[Mi]{miller}
L. Miller: \textit{Curves with invertible Hasse-Witt matrix},
Math. Annalen 197 (1972), 123--127.

\bibitem[Mu]{mumford}
D. Mumford: \textit{Bi-extensions of formal groups}, in Algebraic
Geometry (Internat. Colloq., Tata Inst. Fund. Res., Bombay, 1968),
Oxford Univ. Press, 307--322.

\bibitem[ON]{oort-norman}
P. Norman and F. Oort: \textit{Moduli of abelian varieties}, Ann. of
Math. (2)  112  (1980), no. 3, 413--439.

\bibitem[S]{serre-chebo}
J-P. Serre: \textit{Quelques applications du th\'eor\`eme de densit\'e
  de Chebotarev}, Publ. Math. I.H.E.S 54 (1981), 123--201.

\bibitem[vdW]{vdw}
B.L. van der Waerden: \textit{Moderne algebra}, vol I, Springer 1935.

\bibitem[W]{waterhouse}
W. Waterhouse: \textit{Abelian Varieties over Finite Fields},
Ann. scient. \'Ec. Norm. Sup. 4\`eme s\'erie, 2 (1969), 521--560.

\end{thebibliography}
\end{document}